\documentclass{article}

\newcommand{\HeaderFootnote}{Accepted to be
presented at the IFAC Workshop on Control Applications of
Optimization -- \emph{CAO'2003} -- to be held in Visegr\'{a}d,
Hungary, 30 June -- 2 July 2003, and to appear in the respective
conference proceedings. The date of this version is January 24,
2003.}

\usepackage{amsmath,amssymb,amsthm}

\newtheorem{theorem}{Theorem}[section]
\newtheorem{corollary}[theorem]{Corollary}

\theoremstyle{remark}
\newtheorem{remark}{Remark}[section]
\theoremstyle{definition}
\newtheorem{definition}{Definition}[section]

\theoremstyle{definition}

\begin{document}

\title{Integrals of Motion for Discrete-Time\\
       Optimal Control Problems\footnote{\HeaderFootnote}}

\author{Delfim F. M. Torres\\
        \texttt{delfim@mat.ua.pt}}

\date{Department of Mathematics\\
      University of Aveiro\\
      3810-193 Aveiro, Portugal\\
      \texttt{http://www.mat.ua.pt/delfim}}

\maketitle

%%%%%%%%%%%%%%%%%%%%%%%%%%%

\begin{abstract}
We obtain a discrete time analog of E.~Noether's theorem
in Optimal Control, asserting that integrals of motion
associated to the discrete time Pontryagin Maximum Principle can
be computed from the quasi-invariance properties of the discrete
time Lagrangian and discrete time control system. As corollaries,
results for first-order and higher-order discrete problems
of the calculus of variations are obtained.
\end{abstract}

%%%%%%%%%%%%%%%%%%%%%%%%%%%

\vspace*{0.5cm}

\noindent \textbf{Keywords:} discrete time optimal control,
discrete time calculus of variations, discrete time mechanics,
discrete time Pontryagin extremals,
quasi-invariance up to difference gauge terms,
discrete version of Noether's theorem.

%%%%%%%%%%%%%%%%%%%%%%%%%%%%%%%%%%%%%%%%%%%%%%%%%

\vspace*{0.3cm}


%%%%%%%%%%%%%%%%%%%%%%%%%%%

\section{Introduction}

Most physical systems encountered in nature exhibit
symmetries: there exists appropriate
infinitesimal-parameter family of transformations
which keep the system invariant. From the well-known
theorem of Emmy Noether \cite{JFM46.0770.01,MR53:10538},
one can discover the integrals of motion from those
invariance transformations. Noether's theorem plays a
fundamental role in modern physics, and is usually formulated
in the context of the calculus of variations: from the invariance
properties of the variational integrals, the integrals of
motion of the respective Euler-Lagrange differential equations,
that is, expressions which are preserved along the extremals, are obtained.
The result is, however, much more than a theorem. It is an universal
principle, which can be formalized in a precise statement, as a theorem,
on very different contexts and, for each such context,
under very different assumptions. Let us consider, for example,
classical mechanics or, more generally, the calculus of variations.
Typically, Noether transformations
are considered to be point-transformations (they are considered to be
functions of coordinates and time), but one can consider
more general transformations depending also on velocities and higher
derivatives \cite{ZBL0964.49001} or within the broader context
of dynamical symmetries \cite{MR81c:58054}. For an example of an
integral of motion which comes from an invariance transformation depending
on velocities, see \cite{MoyoLeach2002}.
In most formulations of Noether's principle, the Noether transformations
keep the integral functional invariant
(\textrm{cf. e.g.} \cite[\S 1.5]{MR2000m:49002}).
It is possible, however, to consider transformations of the problem
up to an exact differential (\textrm{cf. e.g.} \cite[p. 73]{MR37:5752}),
called a gauge-term \cite{MR83c:70020}. Once strictly-invariance of the
integral functional is no more imposed, one can think considering additional
terms in the variation of the Lagrangian -- see the quasi-invariance
and semi-invariance notions introduced by the author respectively in
\cite{torresMED2002} and \cite{torresControlo2002}. Formulations
of Noether's principle are possible for problems of the calculus
of variations: on Euclidean spaces (\textrm{cf. e.g.} \cite{MR58:18024})
or on manifolds (\textrm{cf. e.g.} \cite{MR57:13703});
with single or multiple integrals (\textrm{cf. e.g.} \cite{MR95b:58049});
with higher-order derivatives (\textrm{cf. e.g.} \cite{MR54:6786});
with holonomic or nonholonomic constraints (\textrm{cf. e.g.}
\cite[Ch. 7]{MR97a:49001}, \cite{MR2000a:37060}); and so on.
Other contexts for which Noether's theorems are available
include supermechanics \cite{MR96g:58011}, control
systems \cite{MR83k:49054,MR83k:93011}, and optimal control
(see \textrm{e.g.} \cite{MR49:5979,MR1806135,delfim3ncnw,delfimEJC}).
For a survey see \cite{MR96i:49037,torresSpecialJMS}.
Here we are interested in providing a formulation
of the Noether's principle in the discrete time setting.
For a description of discrete time mechanics,
discrete time calculus of variations, and
discrete optimal control see, \textrm{e.g.},
\cite{MR98k:81076,MR98k:81077,MR99m:81275a,MR99m:81275b},
\cite{ZBL0193.07601}, and \cite{MR2001k:49058}.
Illustrative examples of real-life problems which can be modeled
in such framework can be found in \cite[Ch.~8]{SegundaEdicaoMR84g:49002}.
Versions of the Noether's principle for the discrete
calculus of variations, and applicable to
discrete analogues of classical mechanics, appeared
earlier in \cite{MR48:6739,MR48:6741,MR82g:70041,MR95i:58098,%
discreteNoetherPhD,MR99d:70004},
motivated by the advances of numerical and computational methods.
There, the discrete analog of Noether's theorem is obtained
from the discrete analog of the Euler-Lagrange equations.
To the best of our knowledge, no Noether type theorem is available
for the discrete time optimal control setting. One such formulation
is our concern here. The result is obtained from the discrete time
version of the Pontryagin maximum principle. As corollaries,
we obtain generalizations of the previous results for first-order
and higher-order discrete problems of the calculus of variations
which are quasi-invariant and not necessarily invariant.

%%%%%%%%%%%%%%%%%%%%%%%%%%%

\section{Discrete-Time Optimal Control}

Without loss of generality (\textrm{cf.} \cite[\S 2]{MR86c:49028}),
we consider the discrete
optimal control problem in Lagrange form. The time $k$
is a discrete variable: $k \in \mathbb{Z}$.
The horizon consists of $N$ periods,
$k = M, M + 1, \ldots, M + N - 1$, where $M$ and $N$ are
fixed integers, instead of a continuous interval. We look
for a finite control sequence $u(k) \in \mathbb{R}^r$,
$k = M, \ldots, M + N - 1$, and the corresponding state
sequence $x(k) \in \mathbb{R}^n$, $k = M, \ldots, M + N$,
which minimizes or maximizes the sum
\begin{equation*}
J\left[x(\cdot),u(\cdot)\right] = \sum_{k = M}^{M+N-1}
L\left(k,x(k),u(k)\right) \, ,
\end{equation*}
subject to the discrete time control system
\begin{equation}
\label{eq:DCS}
x(k+1) = \varphi\left(k,x(k),u(k)\right) \, , \quad
k = M, \ldots, M + N - 1 \, ,
\end{equation}
the boundary conditions
\begin{equation}
\label{eq:BC}
x(M) = x_{M} \, , \quad x(M+N) = x_{M+N} \, ,
\end{equation}
and the control constraint
\begin{equation*}
u(k) \in \Omega \subseteq \mathbb{R}^r \, , \quad
k = M, \ldots, M + N - 1 \, .
\end{equation*}
A sequence-pair $\left(x(k),u(k)\right)$, $k = M, \ldots, M + N - 1$,
satisfying the recurrence relation \eqref{eq:DCS} and conditions \eqref{eq:BC},
is said to be admissible:
$x(k)$ is an admissible state sequence and $u(k)$ an admissible
control sequence. Functions $L(k,x,u) : \left\{M, \ldots, M + N - 1\right\}
\times \mathbb{R}^n \times \mathbb{R}^r \rightarrow \mathbb{R}$
and $\varphi(k,x,u) : \left\{M, \ldots, M + N - 1\right\}
\times \mathbb{R}^n \times \mathbb{R}^r \rightarrow \mathbb{R}^n$
are assumed to be continuously differentiable with respect to $x$
and $u$ for all fixed $k = M, \ldots, M + N - 1$, and
convex in $u$ for any fixed $k$ and $x$. They are in general nonlinear.
The control constraint set $\Omega$ is assumed to be convex.
The problem is denoted by $(P)$.

\begin{remark}
For continuous optimal control problems,
the convexity assumptions we are imposing are not needed in order
to derive the Pontryagin maximum principle \cite{MR29:3316b}.
This differs from the discrete time optimal control setting.
Our hypothesis can be, however, weakened to directional convexity
or even more weak conditions (see \cite{MR86c:49028},
\cite[\S 8.3]{SegundaEdicaoMR84g:49002} and references in
\cite[Ch. 6]{ZBL0388.49002} and \cite{MR86c:49028}).
\end{remark}

\begin{remark}
It is possible to formulate
problem $(P)$ with the first-order difference equations \eqref{eq:DCS} in terms
of the forward or backward difference operators $\Delta$ or $\nabla$, defined by
$\Delta x(k) = x(k+1) - x(k)$, $\nabla x(k) = x(k) - x(k-1)$.
The results of the paper are written in those terms
in a straightforward way.
\end{remark}

The following theorem provide a first-order necessary
optimality condition (\textrm{cf. e.g.} \cite[\S 3.3.3]{BertsekasVol1},
\cite{MR81i:90157}, \cite[Ch. 6]{ZBL0388.49002}) in the form
of Pontryagin's maximum principle \cite{MR29:3316b}. For a good survey
on the history of the development of maximum principle to the optimization
of discrete time systems, we refer the reader to \cite{MR86c:49028}.

\begin{theorem}[Discrete-Time Maximum Principle]
\label{th:DTMP}
If $\left(x(k),u(k)\right)$ is a minimizer or a maximizer of the problem $(P)$,
then there exists a nonzero sequence-pair $(\psi_0,\psi(k))$,
$k = M+1, \ldots, M+N$, where $\psi_0$ is a constant less or equal
than zero and $\psi(k) \in \mathbb{R}^n$, such that the
sequence-quadruple
\begin{equation*}
\left(x(k),u(k),\psi_0,\psi(k+1)\right) \, , \quad
k = M, \ldots, M+N-1 \, ,
\end{equation*}
satisfies:
\begin{description}
\item[(i)] the Hamiltonian system
\begin{equation}
\label{eq:HamSyst}
\begin{cases}
x(k+1) = \frac{\partial H}{\partial \psi}\left(k,x(k),u(k),\psi_0,\psi(k+1)\right) \, ,
& k = M, \ldots, M+N-1 \, ,\\
\psi(k) = \frac{\partial H}{\partial x}\left(k,x(k),u(k),\psi_0,\psi(k+1)\right) \, ,
& k = M+1, \ldots, M+N-1 \, ;
\end{cases}
\end{equation}
\item[(ii)] the maximality condition
\begin{equation}
\label{eq:MaxCond}
H\left(k,x(k),u(k),\psi_0,\psi(k+1)\right) = \max_{u \in \Omega}
H\left(k,x(k),u,\psi_0,\psi(k+1)\right) \, ,
\end{equation}
$k = M, \ldots, M+N-1$;
\end{description}
with the Hamiltonian
\begin{equation*}
H\left(k,x,u,\psi_0,\psi\right) = \psi_0 L(k,x,u) + \psi \cdot \varphi(k,x,u) \, .
\end{equation*}
\end{theorem}

\begin{remark}
The first equation in the Hamiltonian system is just
the control system \eqref{eq:DCS}. The second equation
in the Hamiltonian system is known as the adjoint system.
The multipliers $\psi(\cdot)$ are called adjoint multipliers
or co-state variables.
\end{remark}

\begin{remark}
In the absence of the initial conditions $x(M) = x_{M}$ and/or terminal
conditions $x(M+N) = x_{M+N}$, there corresponds additional conditions in the
Discrete-Time Maximum Principle called transversality conditions.
Our version of Noether's theorem only require the use
of the adjoint system and maximality condition. Therefore,
the result is valid under all types of boundary conditions
under consideration.
\end{remark}

\begin{definition}
A sequence-quadruple $\left(x(k),u(k),\psi_0,\psi(k+1)\right)$,
$k = M,\ldots,M+N-1$, $\psi_0 \le 0$, satisfying the
Hamiltonian system and the maximality condition, is called
an extremal for problem $(P)$. An extremal is said to be
normal if $\psi_0 \ne 0$ and abnormal if $\psi_0 = 0$.
\end{definition}

\begin{remark}
As we will see on Section~\ref{sec:DCV}, there are no abnormal extremals
both for first-order and higher-order discrete problems of the calculus
of variations. In particular, there are no abnormal extremals for problems
of discrete time mechanics. For our general problem $(P)$, however,
abnormal extremals do exist. In fact, they happen to occur frequently.
For a throughout study of abnormal extremals see \cite{MR1845332}.
\end{remark}

%%%%%%%%%%%%%%%%%%%%%%%%%%%

\section{Integrals of Motion}

We obtain a systematic procedure to establish integrals
of motion, \textrm{i.e.}, to establish expressions which
are preserved on the extremals of the discrete optimal
control problem $(P)$, from the (quasi-)invariance properties
of the discrete Lagrangian $L\left(k,x(k),u(k)\right)$
and discrete control system
$x(k+1) = \varphi\left(k,x(k),u(k)\right)$.

\begin{definition}
\label{def:INV}
Let $X : \left\{M,\ldots,M+N-1\right\} \times \mathbb{R}^n
\times \Omega \times \mathcal{B}(0;\varepsilon) \rightarrow \mathbb{R}^n$,
$\varepsilon > 0$,
$\mathcal{B}(0;\varepsilon) =
\left\{s=(s_1,\ldots,s_\rho)| \left\|s\right\| = \sqrt{\sum_{i=1}^{\rho} (s_i)^2}
< \varepsilon \right\}$,
be an infinitesimal $\rho$-parameter transformation such
that for each $k$, $k = M,\ldots,M+N-1$, $X(k,\cdot,\cdot,\cdot)$ is continuously
differentiable with respect to all arguments, and such that
$X(k,x,u,0) = x$ for all $k = M,\ldots,M+N-1$, $x \in \mathbb{R}^n$,
and $u \in \Omega$. If there exists a real function
$\Phi\left(k,x,u,s\right)$ and for all $s \in \mathcal{B}(0;\varepsilon)$
and admissible $\left(x(k),u(k)\right)$ there exists a control sequence
$u(k,s)$, $u(k,0) = u(k)$, such that:
\begin{multline}
\label{eq:INVi}
L\left(k,x(k),u(k)\right) + \Delta \Phi\left(k,x(k),u(k),s\right)
+ \delta\left(k,x(k),u(k),s\right) \\
= L\left(k,X\left(k,x(k),u(k),s\right),u(k,s)\right) \, ,
\end{multline}
\begin{multline}
\label{eq:INVii}
X\left(k+1,x(k+1),u(k+1),s\right) + \delta\left(k,x(k),u(k),s\right) \\
= \varphi\left(k,X\left(k,x(k),u(k),s\right),u(k,s)\right) \, ,
\end{multline}
for each $k = M,\ldots,M+N-1$ and where $\delta(k,x,u,s)$ is an
arbitrary function satisfying
\begin{equation}
\label{eq:DdeltasEQ0}
\left.\frac{\partial \delta(k,x,u,s)}{\partial s_i}\right|_{s = 0} = 0 \, ,
\quad i = 1,\ldots,\rho \, ,
\end{equation}
for each $k$, $x$, $u$,
then the problem $(P)$ is said
to be quasi-invariant with respect to the transformation $X(k,x,u,s)$
up to the difference gauge term $\Phi\left(k,x(k),u(k),s\right)$.
\end{definition}

\begin{remark}
In the relation \eqref{eq:INVi},
$\Delta$ is the forward difference operator:
\begin{equation*}
\Delta \Phi\left(k,x(k),u(k),s\right) =
\Phi\left(k+1,x(k+1),u(k+1),s\right) - \Phi\left(k,x(k),u(k),s\right) \, .
\end{equation*}
\end{remark}

\begin{remark}
When $\delta \equiv 0$ and $\Phi \equiv 0$, we have (strict-)invariance.
The term \emph{quasi-invariant} refers to the possibility
of $\delta$ to be different from zero.
\end{remark}

\begin{theorem}[Discrete-Time Noether Theorem]
\label{Th:MainResult}
If $(P)$ is quasi-invariant with respect to the $\rho$-parameter transformation
$X$ up to the difference gauge term $\Phi$, in the sense of
Definition~\ref{def:INV}, then all its extremals
$\left(x(k),u(k),\psi_0,\psi(k)\right)$,
$k = M,\ldots,M+N-1$, satisfy the following $\rho$ expressions
($i = 1,\ldots,\rho$):
\begin{equation*}
\psi_0 \left.\frac{\partial}{\partial s_i}
\Phi\left(k,x(k),u(k),s\right)\right|_{s = 0}
+ \psi(k) \cdot \left.\frac{\partial}{\partial s_i}
X\left(k,x(k),u(k),s\right)\right|_{s = 0}
= \text{constant} \, .
\end{equation*}
\end{theorem}

\begin{remark}
The integrals of motion obtained by Theorem~\ref{Th:MainResult}
are ``momentum'' integrals. Due to the fact that time $k$ is discrete,
one can not vary $k$ continuously and, for that reason, one can not obtain
the ``energy'' integrals as in the continuous optimal control case
(\textrm{cf.} \cite{delfim3ncnw,delfimEJC}).
To address the problem another method needs to be developed.
This will be addressed in a forthcoming paper.
\end{remark}

\begin{remark}
Together with the continuous results in \cite{delfim3ncnw,delfimEJC},
Theorem~\ref{Th:MainResult} provides a framework to obtain
a generalization of Noether's theorem for hybrid-systems.
This and related questions are under study and will be addressed elsewhere.
\end{remark}

\begin{proof}
Let $\left(x(k),u(k),\psi_0,\psi(k)\right)$ be an extremal
for problem $(P)$. Differentiating \eqref{eq:INVi} and
\eqref{eq:INVii} with respect to the parameter $s_i$,
$i = 1,\ldots,\rho$, and setting $s = \left(s_1,\ldots,s_{\rho}\right) = 0$,
we get (recall \eqref{eq:DdeltasEQ0} and that $X\left(k,x(k),u(k),0\right)
= x(k)$, $u(k,0) = u(k)$):
\begin{multline}
\label{eq:DsL}
\Delta \left.\frac{\partial}{\partial s_i}
\Phi\left(k,x(k),u(k),s\right)\right|_{s=0}
= \frac{\partial L}{\partial x}\left(k,x(k),u(k)\right)
\cdot \left.\frac{\partial}{\partial s_i} X\left(k,x(k),u(k),s\right)\right|_{s=0} \\
+ \frac{\partial L}{\partial u}\left(k,x(k),u(k)\right)
\cdot \left.\frac{\partial}{\partial s_i} u\left(k,s\right)\right|_{s=0} \, ,
\end{multline}
\begin{multline}
\label{eq:DsFi}
\left.\frac{\partial}{\partial s_i}
X\left(k+1,x(k+1),u(k+1),s\right)\right|_{s=0} \\
= \frac{\partial \varphi}{\partial x}\left(k,x(k),u(k)\right)
\cdot \left.\frac{\partial}{\partial s_i} X\left(k,x(k),u(k),s\right)\right|_{s=0} \\
+ \frac{\partial \varphi}{\partial u}\left(k,x(k),u(k)\right)
\cdot \left.\frac{\partial}{\partial s_i} u\left(k,s\right)\right|_{s=0} \, .
\end{multline}
From the adjoint system
$\psi(k) = \frac{\partial H}{\partial x}\left(k,x(k),u(k),\psi_0,\psi(k+1)\right)$,
we know that
\begin{equation*}
-\psi_0 \frac{\partial L}{\partial x}\left(k,x(k),u(k)\right)
= \psi(k+1) \cdot \frac{\partial \varphi}{\partial x}\left(k,x(k),u(k)\right)
- \psi(k) \, ,
\end{equation*}
and multiplying \eqref{eq:DsL} by $-\psi_0$ one obtains:
\begin{multline}
\label{eq:DepMultpsi0}
\psi_0 \left(\Delta \left.\frac{\partial}{\partial s_i}
\Phi\left(k,x(k),u(k),s\right)\right|_{s=0}
- \frac{\partial L}{\partial u}\left(k,x(k),u(k)\right)
\cdot \left.\frac{\partial}{\partial s_i} u\left(k,s\right)\right|_{s=0}\right) \\
+ \left(\psi(k+1) \cdot \frac{\partial \varphi}{\partial x}\left(k,x(k),u(k)\right)
- \psi(k)\right)
\cdot \left.\frac{\partial}{\partial s_i} X\left(k,x(k),u(k),s\right)\right|_{s=0}
= 0 \, .
\end{multline}
As far as $u(k,0) = u(k)$, according to the maximality condition of the
Discrete-Time Maximum Principle the function
\begin{equation*}
s \mapsto \psi_0 L\left(k,x(k),u(k,s)\right)
+ \psi(k+1) \cdot \varphi\left(k,x(k),u(k,s)\right)
\end{equation*}
attains its maximum for $s = 0$. Therefore,
\begin{equation*}
\left.\frac{\partial}{\partial s_i} \left( \psi_0 L\left(k,x(k),u(k,s)\right)
+ \psi(k+1) \cdot \varphi\left(k,x(k),u(k,s)\right)\right)\right|_{s = 0} = 0 \, ,
\end{equation*}
that is,
\begin{multline}
\label{eq:DaCondMax}
\psi_0 \frac{\partial L}{\partial u}\left(k,x(k),u(k)\right)
\cdot \left.\frac{\partial}{\partial s_i} u(k,s)\right|_{s = 0} \\
+ \psi(k+1) \cdot \frac{\partial \varphi}{\partial u}\left(k,x(k),u(k)\right)
\cdot \left.\frac{\partial}{\partial s_i} u(k,s)\right|_{s = 0} = 0 \, .
\end{multline}
From \eqref{eq:DepMultpsi0} and \eqref{eq:DaCondMax} it comes
\begin{multline*}
\psi_0 \Delta \left.\frac{\partial}{\partial s_i}
\Phi\left(k,x(k),u(k),s\right)\right|_{s=0} \\
+ \left(\psi(k+1) \cdot \frac{\partial \varphi}{\partial x}\left(k,x(k),u(k)\right)
 - \psi(k)\right) \cdot \left.\frac{\partial}{\partial s_i}
 X\left(k,x(k),u(k),s\right)\right|_{s=0} \\
+ \psi(k+1) \cdot \frac{\partial \varphi}{\partial u}\left(k,x(k),u(k)\right)
\cdot \left.\frac{\partial}{\partial s_i} u(k,s)\right|_{s = 0} = 0 \, .
\end{multline*}
Using \eqref{eq:DsFi}, this last equality is equivalent to
\begin{equation*}
\Delta\left(\psi_0 \left.\frac{\partial}{\partial s_i}
\Phi\left(k,x(k),u(k),s\right)\right|_{s = 0}
+ \psi(k) \cdot \left.\frac{\partial}{\partial s_i}
X\left(k,x(k),u(k),s\right)\right|_{s = 0}\right) = 0 \, .
\end{equation*}
The proof is complete.
\end{proof}

%%%%%%%%%%%%%%%%%%%%%%%%%%%

\section{An Example}

We now illustrate the use of Theorem~\ref{Th:MainResult} by
the following example ($n = 3$, $r = 2$, $\Omega = \mathbb{R}^2$):
\begin{gather*}
\sum_{k} \left(u_1(k)\right)^2 - \left(u_2(k)\right)^2 \longrightarrow \textrm{extr} \, , \\
\begin{cases}
x_1(k+1) = x_2(k) + u_1(k) \, , \\
x_2(k+1) = x_1(k) + u_2(k) \, , \\
x_3(k+1) = x_2(k) u_1(k) \, ,
\end{cases}
\end{gather*}
subject to fixed endpoints.
In this case the Hamiltonian is given by
\begin{multline*}
H(x_1,x_2,u_1,u_2,\psi_0,\psi_1,\psi_2,\psi_3) \\
= \psi_0 L(u_1,u_2) + \psi_1 \varphi_1(x_2,u_1)
+ \psi_2 \varphi_2(x_1,u_2) + \psi_3 \varphi_3(x_2,u_1) \, ,
\end{multline*}
with $L(u_1,u_2) = u_1^2 - u_2^2$, $\varphi_1(x_2,u_1) = x_2 + u_1$,
$\varphi_2(x_1,u_2) = x_1 + u_2$, and $\varphi_3(x_2,u_1) = x_2 u_1$.
From the adjoint system we get the evolution equations
\begin{equation*}
\begin{split}
\psi_1(k) &= \psi_2(k+1) \, ,\\
\psi_2(k) &= \psi_1(k+1) + \psi_3(k+1) u_1(k) \, , \\
\psi_3(k) &= 0 \, ,
\end{split}
\end{equation*}
while from the maximality conditions we get ($\psi_3 = 0$)
\begin{equation*}
\begin{split}
\psi_1(k+1) &= -2 \psi_0 u_1(k) \, ,\\
\psi_2(k+1) &= 2 \psi_0 u_2(k) \, .
\end{split}
\end{equation*}
There are no abnormal extremals for the problem,
and one can fix $\psi_0 = -\frac{1}{2}$. The extremals
are obtained solving five difference-equations of order one,
\begin{equation*}
\begin{cases}
x_1(k+1) = x_2(k) + \psi_1(k+1) \, , \\
x_2(k+1) = x_1(k) - \psi_2(k+1) \, , \\
x_3(k+1) = x_2(k) \psi_1(k+1) \, , \\
\psi_1(k+1) = \psi_2(k) \, , \\
\psi_2(k+1) = \psi_1(k) \, ,
\end{cases}
\end{equation*}
together with the boundary conditions (or the transversality conditions),
by standard techniques.
On the other hand, the problem is quasi-invariant with respect to
the one-parameter ($\rho = 1$) transformations
\begin{equation*}
\begin{split}
X_1\left(x_1(k),s\right) &= x_1(k)+2s \, , \\
X_2\left(x_2(k),s\right) &= x_2(k)+s \, , \\
X_3\left(x_1(k),x_3(k),s\right) &= x_3(k) + s x_1(k) \, ,
\end{split}
\end{equation*}
up to the difference gauge term
$\Phi\left(x_1(k),x_2(k),s\right) = 2\left(x_1(k)+x_2(k)\right) s$.
To see that we choose
\begin{equation*}
u_1(k,s) = u_1(k) + s \, , \quad
u_2(k,s) = u_2(k) - s \, ,
\end{equation*}
in the Definition~\ref{def:INV}. We notice that
$X_1\left(x_1(k),0\right) = x_1(k)$,
$X_2\left(x_2(k),0\right) = x_2(k)$,
$X_3\left(x_1(k),x_3(k),0\right) = x_3(k)$,
$u_1(k,0) = u_1(k)$, and $u_2(k,0) = u_2(k)$.
Direct verifications show that the quasi-invariance
conditions are satisfied:
\begin{equation*}
\begin{split}
L(u_1(k,s),&u_2(k,s))
= \left(u_1(k)\right)^2 - \left(u_2(k)\right)^2 + 2\left(u_1(k)+u_2(k)\right) s \\
&= L\left(u_1(k),u_2(k)\right) + 2 \left(x_1(k+1)-x_2(k)+x_2(k+1)-x_1(k)\right) s \\
&= L\left(u_1(k),u_2(k)\right) + \Delta \Phi\left(x_1(k),x_2(k),s\right) \, ,
\end{split}
\end{equation*}
\begin{gather*}
\begin{split}
\varphi_1\left(X_2\left(x_2(k),s\right),u_1(k,s)\right)
&= x_2(k) + u_1(k) + 2s = x_1(k+1) + 2s \\
&= X_1\left(x_1(k+1),s\right) \, , \\
\end{split}\\
\begin{split}
\varphi_2\left(X_1\left(x_1(k),s\right),u_2(k,s)\right)
&= x_1(k) + u_2(k) + s = x_2(k+1) + s \\
&= X_2\left(x_2(k+1),s\right) \, , \\
\end{split}\\
\begin{split}
\varphi_3\left(X_2\left(x_2(k),s\right),u_1(k,s)\right)
&= \left(x_2(k)+s\right) \left(u_1(k)+s\right) \\
&= x_2(k) u_1(k)
+ s \left(x_2(k) + u_1(k)\right) + s^2 \\
&= x_3(k+1) + s x_1(k+1) + \delta(s) = X_3(k+1) + \delta(s) \, .
\end{split}
\end{gather*}
By Theorem~\ref{Th:MainResult} we obtain the following conservation
law for the problem:
\begin{equation}
\label{eq:FIMEx}
2 \psi_0 \left(x_1(k)+x_2(k)\right) + 2 \psi_1(k)
+ \psi_2(k) + \psi_3(k) x_1(k) = \text{constant} \, .
\end{equation}
Using the information from the discrete time maximum principle,
condition \eqref{eq:FIMEx} is equivalent to
\begin{equation}
\label{eq:FIMExCondNec}
\left(x_1(k)+x_2(k)\right) + 2 u_2(k) - u_1(k) = \text{constant} \, .
\end{equation}
The conservation law \eqref{eq:FIMExCondNec} is a necessary
optimality condition. It is trivially satisfied choosing the
control variables according to:
\begin{gather*}
u_1(k) = x_1(k) \, ,\\
u_2(k) = - \frac{1}{2} x_2(k) \, .
\end{gather*}
An extremal is then obtained with the co-state variables given by
\begin{gather*}
\psi_1(k) = - u_2(k) = \frac{1}{2} x_2(k) \, , \\
\psi_2(k) = u_1(k) = x_1(k) \, , \\
\psi_3(k) = 0 \, .
\end{gather*}

%%%%%%%%%%%%%%%%%%%%%%%%%%%

\section{Discrete Calculus of Variations}
\label{sec:DCV}

We now obtain a discrete Noether's theorem for the problems
of the discrete time calculus of variations
which are quasi-invariant with respect to
infinitesimal transformations having $\rho$ parameters, $\rho \ge 1$,
up to a difference gauge term.

%%%%%%%%%%%%%

\subsection{The fundamental Problem}

The fundamental problem in the discrete calculus
of variations is a special case of our problem $(P)$:
$r = n$; no restrictions on the controls
($\Omega = \mathbb{R}^n$); $\varphi(k,x,u) = u$.
The problem is then to determine a finite sequence
$x(k) \in \mathbb{R}^n$, $k = M,\ldots,M+N$,
$x(M) = x_{M}$, $x(M+N) = x_{M+N}$,
for which the discrete cost function
\begin{equation*}
J\left[x(\cdot)\right] = \sum_{k = M}^{M+N-1}
L\left(k,x(k),x(k+1)\right)
\end{equation*}
is extremized. The maximality condition in
the Theorem~\ref{th:DTMP} implies in this case the conditions
\begin{equation*}
\frac{\partial H}{\partial u}\left(k,x(k),u(k),\psi_0,\psi(k+1)\right) = 0 \, ,
\quad k = M,\ldots,M+N-1 \, ,
\end{equation*}
that is,
\begin{equation}
\label{eq:NoAbnPB}
\psi(k+1) = -\psi_0 \frac{\partial L}{\partial u}\left(k,x(k),x(k+1)\right) \, ,
\end{equation}
while from the adjoint system one gets
\begin{equation}
\label{eq:SAPBCV}
\psi(k) = \psi_0 \frac{\partial L}{\partial x}\left(k,x(k),x(k+1)\right) \, .
\end{equation}
We note that no abnormal extremals exist for the fundamental
problem of the discrete calculus of variations: $\psi_0 = 0$ implies
that $\psi(k+1)$ is zero for all $k = M,\ldots,M+N-1$, a possibility
excluded by the discrete time maximum principle. So it must be the
case that $\psi_0 \ne 0$. From \eqref{eq:NoAbnPB} and \eqref{eq:SAPBCV},
a necessary condition to have an extremum is that $x(k)$,
$k = M,\ldots,M+N-2$, must satisfy the second-order difference equation
\begin{equation}
\label{eq:DiscELeq}
\frac{\partial L}{\partial x}\left(k+1,x(k+1),x(k+2)\right)
+ \frac{\partial L}{\partial u}\left(k,x(k),x(k+1)\right) = 0 \, .
\end{equation}
Equations \eqref{eq:DiscELeq} share resemblances with the
continuous Euler-Lagrange equations, and are called the
discrete Euler-Lagrange equations.

\begin{definition}
\label{def:INVPBCV}
The discrete Lagrangian $L\left(k,x(k),x(k+1)\right)$
is said to be quasi-invariant with respect to the infinitesimal $\rho$-parameter
transformation $X(k,x,u,s)$, $s = \left(s_1,\ldots,s_{\rho}\right)$,
$\left\|s\right\| < \varepsilon$,
$X(k,x,u,0) = x$ for all $k = M,\ldots,M+N-1$, $x, u \in \mathbb{R}^n$,
up to the difference gauge term $\Phi\left(k,x(k),x(k+1),s\right)$,
if for each $k = M,\ldots,M+N-2$
\begin{multline}
\label{eq:INViPBCV}
L\left(k,x(k),x(k+1)\right) + \Delta \Phi\left(k,x(k),x(k+1),s\right)
+ \delta\left(k,x(k),x(k+1),s\right)\\
= L\left(k,X\left(k,x(k),x(k+1),s\right),X\left(k+1,x(k+1),x(k+2),s\right)\right) \, ,
\end{multline}
where $\delta(\cdot,\cdot,\cdot,\cdot)$ is a function satisfying \eqref{eq:DdeltasEQ0}.
\end{definition}

\begin{corollary}
\label{CorMainResultPBCV}
If $L\left(k,x(k),x(k+1)\right)$ is quasi-invariant with respect
to the $\rho$-parameter transformation
$X$ up to the difference gauge term $\Phi$, in the sense of
Definition~\ref{def:INVPBCV}, then all solutions $x(k)$,
$k = M,\ldots,M+N-2$, of the discrete Euler-Lagrange difference
equation \eqref{eq:DiscELeq} satisfy
\begin{multline*}
\frac{\partial L}{\partial u}\left(k-1,x(k-1),x(k)\right)
\cdot \left.\frac{\partial}{\partial s_i}
X\left(k,x(k),x(k+1),s\right)\right|_{s = 0} \\
- \left.\frac{\partial}{\partial s_i}
\Phi\left(k,x(k),x(k+1),s\right)\right|_{s = 0} = \text{constant} \, ,
\end{multline*}
$i = 1,\ldots,\rho$.
\end{corollary}

%%%%%%%%%%%%%

\subsection{Higher Order Discrete Problems}

Let us now consider the problem of optimizing
\begin{equation}
\label{eq:Pr:HODP}
\sum_{k} L\left(k,x(k),x(k+1),\ldots,x(k+m)\right) \, ,
\end{equation}
where $L(k,x^0,x^1,\ldots,x^{m})$ is continuously differentiable
with respect to all variables.
This problem is analogous to the continuous problems
of the calculus of variations for which the Lagrangian
$L$ depends on higher-order derivatives. It is easily written
in the optimal control form $(P)$. Introducing the notation
\begin{gather*}
x^0(k) = x(k) \, , \\
x^1(k) = x(k+1) \, , \\
\vdots \\
x^{m-1}(k) = x(k+m-1) \, , \\
u(k) = x(k+m) \, ,
\end{gather*}
one gets:
\begin{gather*}
\sum_{k} L\left(k,x^0(k),\ldots,x^{m-1}(k),u(k)\right)
\longrightarrow \textrm{extr} \, , \\
\begin{cases}
x^0(k+1) = x^1(k) \, , \\
x^1(k+1) = x^2(k) \, , \\
\vdots \\
x^{m-2}(k+1) = x^{m-1}(k) \, , \\
x^{m-1}(k+1) = u(k) \, .
\end{cases}
\end{gather*}
The Hamiltonian is given by
\begin{equation*}
H = \psi_0 L\left(k,x^0,\ldots,x^{m-1},u\right)
+ \left(\sum_{j=0}^{m-2} \psi^j \cdot x^{j+1}\right) + \psi^{m-1} u \, .
\end{equation*}
From the maximality condition
\begin{equation}
\label{eq:MCHighOrderPr}
\psi^{m-1}(k+1) = - \psi_0
\frac{\partial L}{\partial x^m}\left(k,x^0(k),\ldots,x^{m-1}(k),u(k)\right) \, ,
\end{equation}
while from the adjoint system
\begin{gather}
\psi^0(k) =  \psi_0
\frac{\partial L}{\partial x^0}\left(k,x^0(k),\ldots,x^{m-1}(k),u(k)\right)
\, , \label{eq:ASHighOrderPr1} \\
\psi^{j}(k) = \psi_0
\frac{\partial L}{\partial x^j}\left(k,x^0(k),\ldots,x^{m-1}(k),u(k)\right)
+ \psi^{j-1}(k+1) \, , \label{eq:ASHighOrderPr2}
\end{gather}
$j = 1,\ldots,m-1$.
From \eqref{eq:MCHighOrderPr}, \eqref{eq:ASHighOrderPr1}, and
\eqref{eq:ASHighOrderPr2}, we conclude that:
similarly to the fundamental problem of the calculus of variations,
no abnormal extremals exist in the higher order case;
the equation
\begin{equation}
\label{eq:EulerPoissonOCnotation}
\sum_{j=0}^{m} \frac{\partial L}{\partial
x^j}\left(k+m-j,x^0(k+m-j),\ldots,x^{m-1}(k+m-j),u(k+m-j)\right) = 0
\end{equation}
holds. Going back to the initial notation,
\eqref{eq:EulerPoissonOCnotation} is nothing more than
the discrete Euler-Poisson equation of order $2m$ for the
$m$-th order discrete problem of the calculus of variations
\eqref{eq:Pr:HODP}:
\begin{equation}
\label{eq:EulerPoisson}
\sum_{j=0}^{m} \frac{\partial L}{\partial
x^j}\left(k+m-j,x(k+m-j),\ldots,x(k+2m-1-j),x(k+2m-j)\right) = 0 \, .
\end{equation}

\begin{definition}
\label{def:INVPHOCV}
The discrete Lagrangian $L\left(k,x(k),\ldots,x(k+m)\right)$
is said to be quasi-invariant with respect to the infinitesimal $\rho$-parameter
transformation $X(k,x^0,\ldots,x^m,s)$, $s = \left(s_1,\ldots,s_{\rho}\right)$,
$\left\|s\right\| < \varepsilon$,
$X(k,x^0,\ldots,x^m,0) = x^0$ for all $k$, and $x^j$, $j = 0,\ldots,m$,
up to the difference gauge term $\Phi\left(k,x(k),\ldots,x(k+m),s\right)$,
if for each $k$
\begin{multline}
\label{eq:INViPHOCV}
L\left(k,x(k),\ldots,x(k+m)\right)
+ \Delta \Phi\left(k,x(k),\ldots,x(k+m),s\right)
+ \delta\left(k,x(k),\ldots,x(k+m),s\right) \\
= L\left(k,X\left(k,x(k),\ldots,x(k+m),s\right),\ldots,
X\left(k+m,x(k+m),\ldots,x(k+2m),s\right)\right) \, ,
\end{multline}
where $\frac{\partial \delta}{\partial s_i} = 0$, $i = 1,\ldots,\rho$.
\end{definition}

\begin{corollary}
\label{CorMainResultPHOCV}
If $L\left(k,x(k),\ldots,x(k+m)\right)$ is quasi-invariant with respect
to the $\rho$-parameter transformation
$X$ up to the difference gauge term $\Phi$, in the sense of
Definition~\ref{def:INVPHOCV}, then all solutions $x(k)$
of the discrete Euler-Poisson difference
equation \eqref{eq:EulerPoisson} satisfy
\begin{multline}
\label{eq:LeiConsPHOCV}
\left.\frac{\partial}{\partial s_i}
\Phi\left(k,x(k),\ldots,x(k+m),s\right)\right|_{s = 0} \\
+ \sum_{j=0}^{m-1} \sum_{l=0}^{j}
\frac{\partial L}{\partial x^l}\left(k+j-l,x(k+j-l),\ldots,x(k+j-l+m)\right) \\
\cdot \left.\frac{\partial}{\partial s_i}
X\left(k+j,x(k+j),\ldots,x(k+j+m),s\right)\right|_{s = 0}
= \text{constant} \, ,
\end{multline}
$i = 1,\ldots,\rho$.
\end{corollary}

In the case $m=1$ the discrete Euler-Poisson equation \eqref{eq:EulerPoisson}
reduces to the discrete Euler-Lagrange equation \eqref{eq:DiscELeq}, and
the conservation law \eqref{eq:LeiConsPHOCV} reduces to
\begin{multline*}
\left.\frac{\partial}{\partial s_i}
\Phi\left(k,x(k),x(k+1),s\right)\right|_{s = 0} \\
+ \frac{\partial L}{\partial x^0}\left(k,x(k),x(k+1)\right)
\cdot \left.\frac{\partial}{\partial s_i}
X\left(k,x(k),x(k+1),s\right)\right|_{s = 0}
= \text{constant} \, ,
\end{multline*}
or, which is the same,
\begin{multline*}
\left.\frac{\partial}{\partial s_i}
\Phi\left(k,x(k),x(k+1),s\right)\right|_{s = 0} \\
- \frac{\partial L}{\partial x^1}\left(k-1,x(k-1),x(k)\right)
\cdot \left.\frac{\partial}{\partial s_i}
X\left(k,x(k),x(k+1),s\right)\right|_{s = 0}
= \text{constant} \, .
\end{multline*}
This is precisely the conservation law given
by Corollary~\ref{CorMainResultPBCV}.

%%%%%%%%%%%%%%%%%%%%%%%%%%%%%%%%%%%%%%%%%%%%%%%%%%%%%%%%%%%%%%%%%%%

%%%%%%%%%%%%%%%%%%%%%%%%%%%%%%%%%%%%%%%%%%%%%%%%%%%%%%%%%%%%%%%%%%%

\end{document}